\documentclass[11pt, oneside]{article}   	% use "amsart" instead of "article" for AMSLaTeX format
\usepackage{geometry}                		% See geometry.pdf to learn the layout options. There are lots.
\usepackage{graphicx}				% Use pdf, png, jpg, or eps§ with pdflatex; use eps in DVI mode
\usepackage{times}
\usepackage{graphicx}
\usepackage{amssymb}
\usepackage{nomencl}
\usepackage{verbatim}
\usepackage{amsmath}
\usepackage{amsfonts}
\usepackage{amssymb}
\usepackage{latexsym}
\usepackage{pdfpages}
\usepackage{longtable}
\usepackage{ragged2e}
\usepackage{array}
\usepackage{textcomp}
\usepackage{soul}
\usepackage[titletoc,title]{appendix}
\usepackage{color}
\usepackage{epstopdf}
\usepackage[margin=10pt,font=small,labelfont=bf]{caption}
\usepackage{subcaption}
\newcommand{\beq}{\begin{eqnarray}}
\newcommand{\eeq}{\end{eqnarray}}

\newcommand{\be}{\begin{enumerate}}
\newcommand{\ee}{\end{enumerate}}
\newcommand{\bi}{\begin{itemize}}
\newcommand{\ei}{\end{itemize}}
\newcommand{\bc}{\begin{center}}
\newcommand{\ec}{\end{center}}
\newcommand{\bt}{\begin{tabbing}}
\newcommand{\et}{\end{tabbing}}

\title{On the shape of air-liquid interfaces with surface tension that\\ bound rigidly-rotating liquids in partially filled containers}
\author{\vspace{-0.5em}Enrique Ram\'{e}\\
\textit{Madrid, 28034 Spain}\\
\vspace{-0.5em}
Steven J. Weinstein\\\vspace{-0.5em}
\textit{Department of Chemical Engineering}\\\vspace{-0.5em}
\textit{School of Mathematical Sciences}\\
\textit{Rochester Institute of Technology, Rochester, NY 14623, USA}\\
\vspace{-0.5em}
Nathaniel S. Barlow\\\vspace{-0.5em}
\textit{School of Mathematical Sciences}\\
\textit{Rochester Institute of Technology, Rochester, NY 14623, USA}
}
\date{}							% Activate to display a given date or no date

\begin{document}
\maketitle

\begin{abstract}
The interface shape of a fluid in rigid body rotation about its axis and partially filling the container is often the subject of a homework problem in the first graduate fluids class. In that problem, surface tension is neglected, the interface shape is parabolic and the contact angle boundary condition is not satisfied in general. When surface tension is accounted for, the shapes exhibit much richer dependencies as a function of rotation velocity. We analyze steady interface shapes in rotating right-circular cylindrical containers under rigid body rotation in zero gravity. We pay especial attention to shapes near criticality, in which the interface, or part thereof, becomes straight and parallel to the axis of rotation at certain specific rotational speeds. We examine geometries where the container is axially infinite and derive properties of their solutions. We then examine in detail two special cases of menisci in a cylindrical container: a meniscus spanning the cross section; and a meniscus forming a bubble. In each case we develop exact solutions for the respective axial lengths as infinite series in powers of appropriate rotation parameters; and we find the respective asymptotic behaviors as the shapes approach their critical configuration. Finally we apply the method of asymptotic approximants to yield analytical expressions for the axial lengths of menisci over the whole range of rotation speeds. In this application, the analytical solution is employed to examine errors introduced by the assumption that the interface is a right circular cylinder; this assumption is key to the spinning bubble method used to measure surface tension.\end{abstract}
\newpage

\section{Introduction}
\begin{figure}[h]
\centering
\begin{subfigure}{0.4\textwidth}
  \centering
  \includegraphics[width=2.5in]{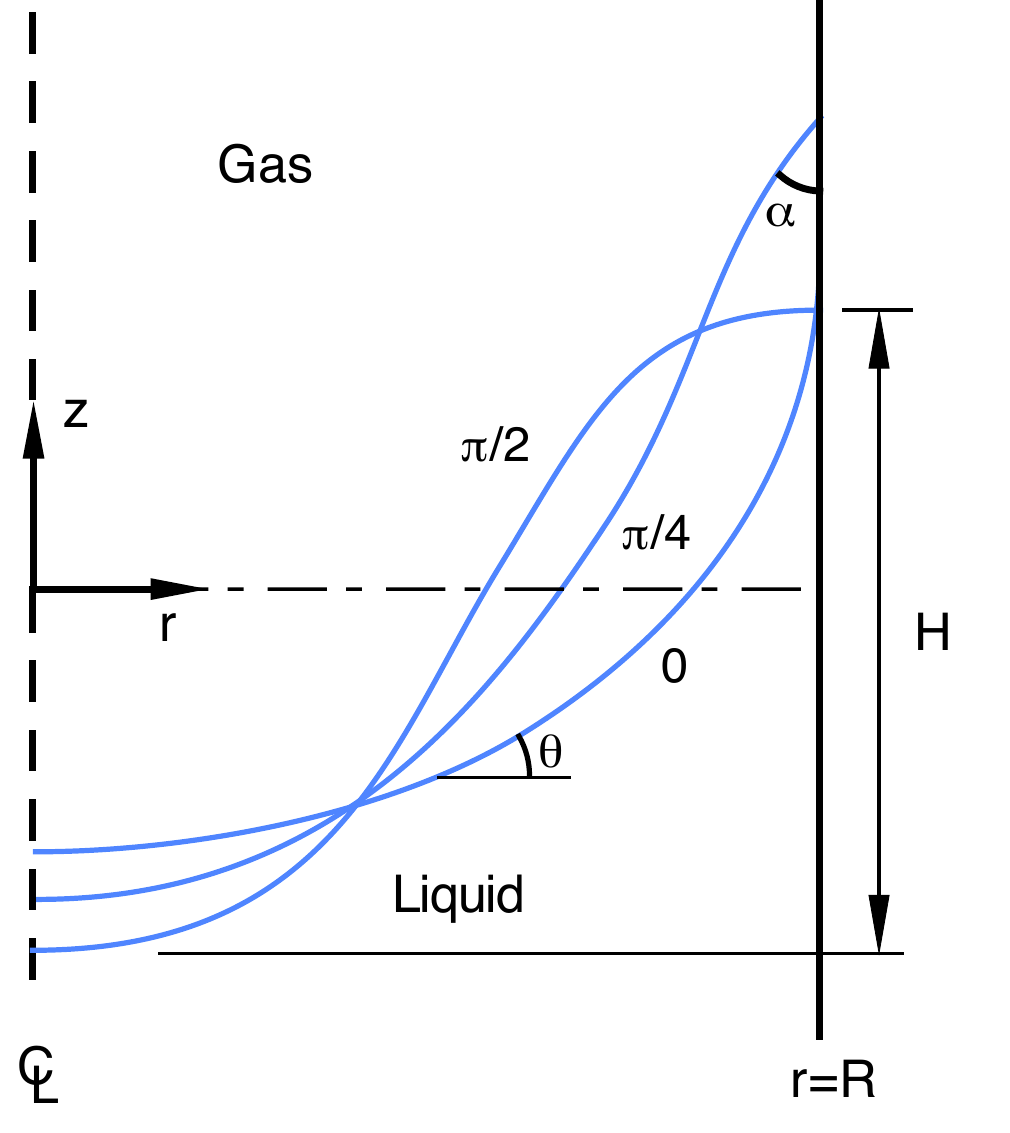}
  \caption{}
  \label{shapesschem}
\end{subfigure}%
\begin{subfigure}{0.3\textwidth}
  \centering
  \includegraphics[width=1.0in]{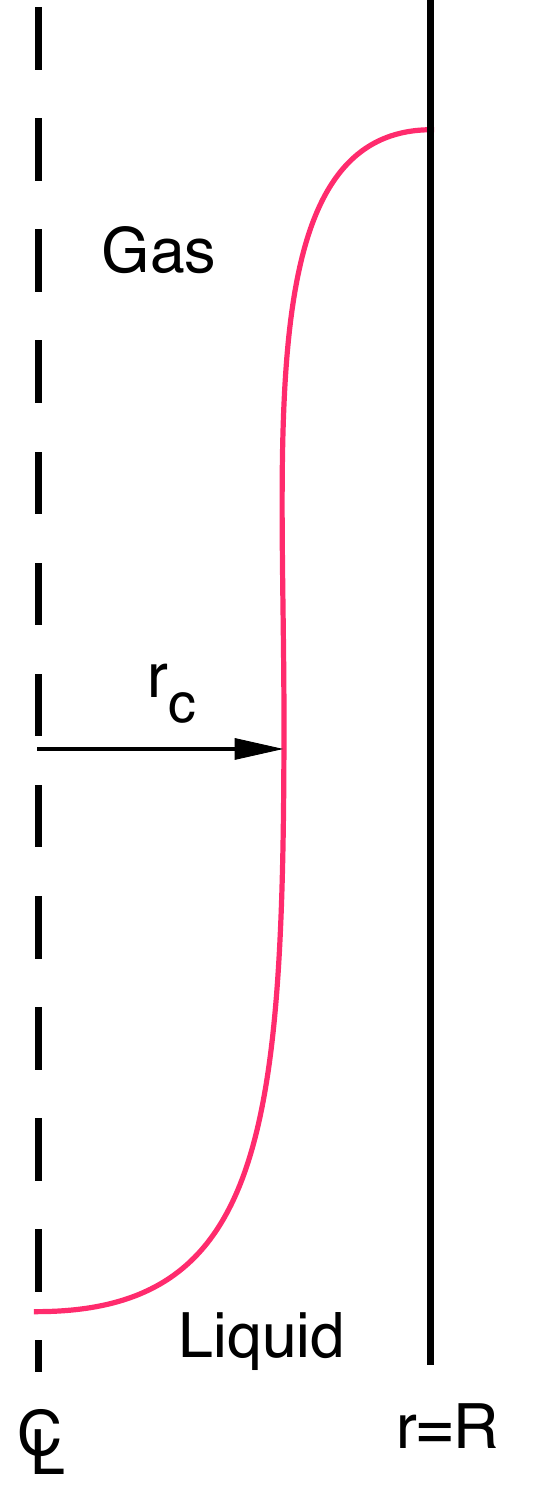}
  \caption{}
  \label{criticalschem}
\end{subfigure}
\begin{subfigure}{0.22\textwidth}
  \centering
  \includegraphics[width=1.61in]{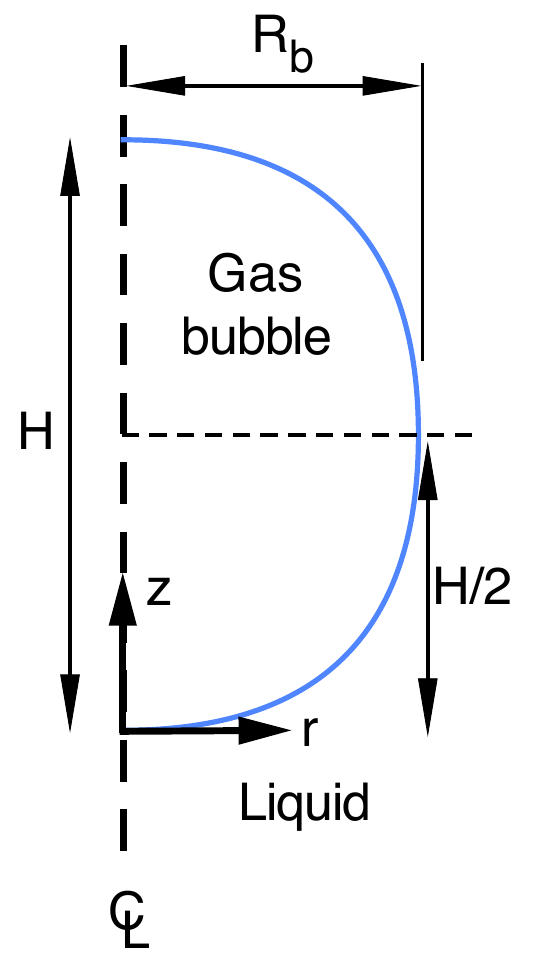}
  \caption{}
  \label{bubbleschem}
\end{subfigure}
\caption{Schematics of surface configurations with fluids in rigid body rotation about the axis of a cylindrical container. Liquid is below the interface. (a): typical shapes at $\lambda<\lambda_c$, for wall contact angles, $\alpha=0, \pi/4$ and $\pi/2$, as marked. $\theta$ is the slope angle of the interface. $z=0$ is chosen so that the volumes of liquid above and below are equal. $H$ is the axial meniscus length. (b): Shape with $\lambda$ close to $\lambda_c$, for 90-degree contact angle. Most of the surface is a straight circular cylinder of radius $r_c$. $H$ continually increases as $\lambda\to\lambda_c$. (c): Spinning bubble. $R_b$ is the maximum radius and $H$ the total axial length. The half-bubble on either side of the equator is mathematically identical to the $\alpha=0$ case in fig. \ref{shapesschem}.}
\label{schematsurf}
\end{figure}

The shapes of fluid interfaces in rigid body rotation have been well studied, with the spinning bubble tensiometer being a notable example, see \cite{vonnegut1942, joseph1994}. Since about the mid 1950s, interest in such problems has grown with the need to engineer fluid containers in zero gravity for spacecraft, where guaranteeing a known location for the liquid phase is crucial for rockets to fire properly. Seebold \cite{seebold1965} performed a stability analysis of menisci in circular cylindrical containers in rigid-body rotation with arbitrary axial gravity and contact angle. He derived stability limits by a variational analysis of the Hamilton principle. Later Preziosi \& Joseph \cite{preziosi1987} analyzed the stability of periodic interface shapes in rigid body rotation by minimization of an energy potential under conditions of negligible gravity, obtaining results that are wholly consistent with those of Seebold \cite{seebold1965}.
\\
\\
Similar to a static meniscus in a gravitational field, the shape of the interface between two immiscible fluids in rigid body rotation with angular velocity $\omega$, density difference $\Delta\rho\equiv \rho_1-\rho_2$ ($>0$) and surface tension $\sigma$, depends on the rotational Bond number $\lambda\equiv\Delta\rho\,\omega^2 d^3/\sigma$, where $d$ is some appropriate length that depends on the geometry being considered. The studies cited above show that a critical value $\lambda_c$ exists such that, when $\lambda\to\lambda_c$ from below, the interface undergoes a critical transition, with an outcome that depends on the configuration of the fluid body and the particular container. Figure \ref{schematsurf} illustrates some typical configurations.
\\
\\
More specifically, Seebold \cite{seebold1965} showed that menisci having finite axial length and spanning the entire cross section of a cylinder with contact angle $\alpha$ at the cylinder wall (such as shown in fig. \ref{shapesschem}), exist in zero gravity under rigid body rotation for $\lambda<4 f(\alpha)$, where $f(\alpha)$ is derived later in eq. (\ref{criticalrlambda}) and satisfies $f(0)=1$. When $\lambda$ approaches $4 f(\alpha)$, the meniscus reorganizes with a divergent axial length into the shape of a straight cylinder (fig. \ref{criticalschem}). Preziosi \& Joseph \cite{preziosi1987} showed that non-straight periodic interfaces that develop  in rigid body motions for $\lambda < 4$ become straight cylinders when when $\lambda=4$. This result is entirely consistent with the analysis of Seebold \cite{seebold1965}, who collected extensive experimental evidence (and also intuited without proof) that shapes become straight cylinders when $\lambda\to\lambda_c$. Ross \cite{ross1968} studied the shapes of rotating drops and bubbles (as shown in fig. \ref{bubbleschem}), and obtained several important results and interpretations that coincide with those of Seebold \cite{seebold1965} and Preziosi \& Joseph \cite{preziosi1987} in the case of bubbles, and with those of Chandrasekhar \cite{chandrasekhar1965} in the case of drops.
\\
\\
In this work we derive new features of menisci in rigid body rotation. Specifically, our goal is to develop analytical tools to describe the axial meniscus length, denoted as $H$, in configurations of practical relevance as depicted in fig. \ref{schematsurf}. We analyze two configurations: 1) a meniscus spanning the cylinder cross section, as in figs. \ref{shapesschem} and \ref{criticalschem}; and 2) a meniscus forming a bubble whose axis coincides with the axis of rotation, as in fig. \ref{bubbleschem}. The latter geometry corresponds to that of the spinning bubble tensiometer \cite{vonnegut1942}. For each geometry, we examine the respective shapes and axial meniscus lengths as a function of $\lambda$. In particular, we identify the asymptotic behaviors of the divergence in axial meniscus length as $\lambda\to\lambda_c$ and find that the divergence law depends on the meniscus configuration. We then develop exact solutions for these lengths as infinite series in powers of $\lambda$. Since these solutions converge poorly near $\lambda_c$, we apply the method of asymptotic approximants \cite{barlow2017} to describe axial meniscus length uniformly over the whole range $0\le\lambda<\lambda_c$ for both configurations. We end with a discussion of implications relevant the measurement of surface tension using a spinning bubble tensiometer.

\section{Analysis}
\subsection{Formulation}
Consider a liquid of density $\rho$ in contact with a gas of negligible density in a cylindrical container of radius $R$, rotating with angular velocity $\omega$ about its axis in rigid body rotation (fig. \ref{schematsurf}). The gas-liquid interface has surface tension $\sigma$, its location is $z=h(r)$, and it obeys the normal component of the dynamic boundary condition with a pressure whose gradient arises solely from centripetal acceleration.\\
\\
Neglecting the dynamics of the gas and using $d$ as a characteristic length scale, the dimensionless governing equation is given as
\beq
\frac{1}{r}\frac{d}{dr}\left(\frac{r h'}{\sqrt{1+h'^2}}\right)
         =-P_0 - \lambda \frac{r^2}{2},
\label{dbc}
\eeq
where $h'(r)=\tan\theta$ is the slope relative to the $r$-axis ($\theta$ is shown explicitly in fig. \ref{shapesschem}), $P_0$ is the pressure difference across the interface at $r=0$ made dimensionless with $\sigma/d$ and $\lambda$ is the rotational Bond number, defined as 
\beq
\lambda\equiv\frac{\rho\omega^2 d^3}{\sigma}.
\label{lambda}
\eeq
%(If the meniscus contacts an interior surface that prevents it from reaching the center line, then $P_0$ is the fictitious pressure difference at $r=0$, obtained by assuming that the meniscus is present there.) 
In general, the characteristic length $d$ depends on the configuration considered, as sketched in fig. \ref{schematsurf}. If the interface spans the entire radius of the container, then $d=R$; for a bubble wholly surrounded by the liquid, $d$ is the maximum radius of the bubble, $R_b$.\\
\\
Since $h'(r)=\tan\theta$, the left-hand side of eq. (\ref{dbc}) may be written as $d(r\sin\theta)/(r\,dr)$. This equation may be integrated once to obtain an expression for $\sin\theta$ as a function of $r$ using two boundary conditions. One boundary condition accounts for the axial symmetry and is imposed to all the cases studied here:
\beq
\theta(0)=0.
\eeq
The other boundary condition is applied at the radial end of the interface. In the case of a meniscus spanning the container, the interface obeys a contact angle condition:
\beq
\theta(1)=\frac{\pi}{2}-\alpha,
\label{canglebc}
\eeq
where $\alpha$ is the contact angle, as depicted in fig. \ref{shapesschem}. In the case of a bubble wholly surrounded by liquid, the second boundary condition is:
\beq
\theta(1)=\frac{\pi}{2}.
\label{bubblebc}
\eeq
Physically, this condition reflects the equatorial symmetry of the bubble; but mathematically it can be seen to be identical to the case of a meniscus with contact angle $\alpha=0$ expressed in eq. (\ref{canglebc}). The shape for arbitrary contact angle $\alpha$ is therefore:
\beq
\frac{h'}{\sqrt{1+h'^2}}=\sin\theta=r\cos\alpha+\frac{\lambda}{8}r(1-r^2),
\label{sint}
\eeq
where, from the preceding discussion, $\alpha$ must be set to zero to describe a bubble. One additional integration determines $h(r)$ after a constraint specific to each configuration is applied. The axial meniscus length $H(\lambda)$ for each interface configuration is computed by integrating $h'(r)$ while taking proper account of the limits of integration. As can be seen in fig. \ref{schematsurf}, $h(r)$ is single-valued in a meniscus spanning the container radius; and multivalued in a spinning bubble; therefore,
\begin{subequations}
\begin{align}
&H(\lambda)=\int_0^1 h'(r) dr:~~~~~~~~\text{meniscus spanning container;}\\
&H(\lambda)=2\int_0^1 h'(r) dr:~~~~~~\text{bubble.}
\end{align}
\label{Hdef}
\end{subequations}

\subsection{Properties of interface shapes satisfying eq. \ref{sint}}\label{properties}
\begin{figure}[t]
\begin{center}
\includegraphics[width=3.5in]{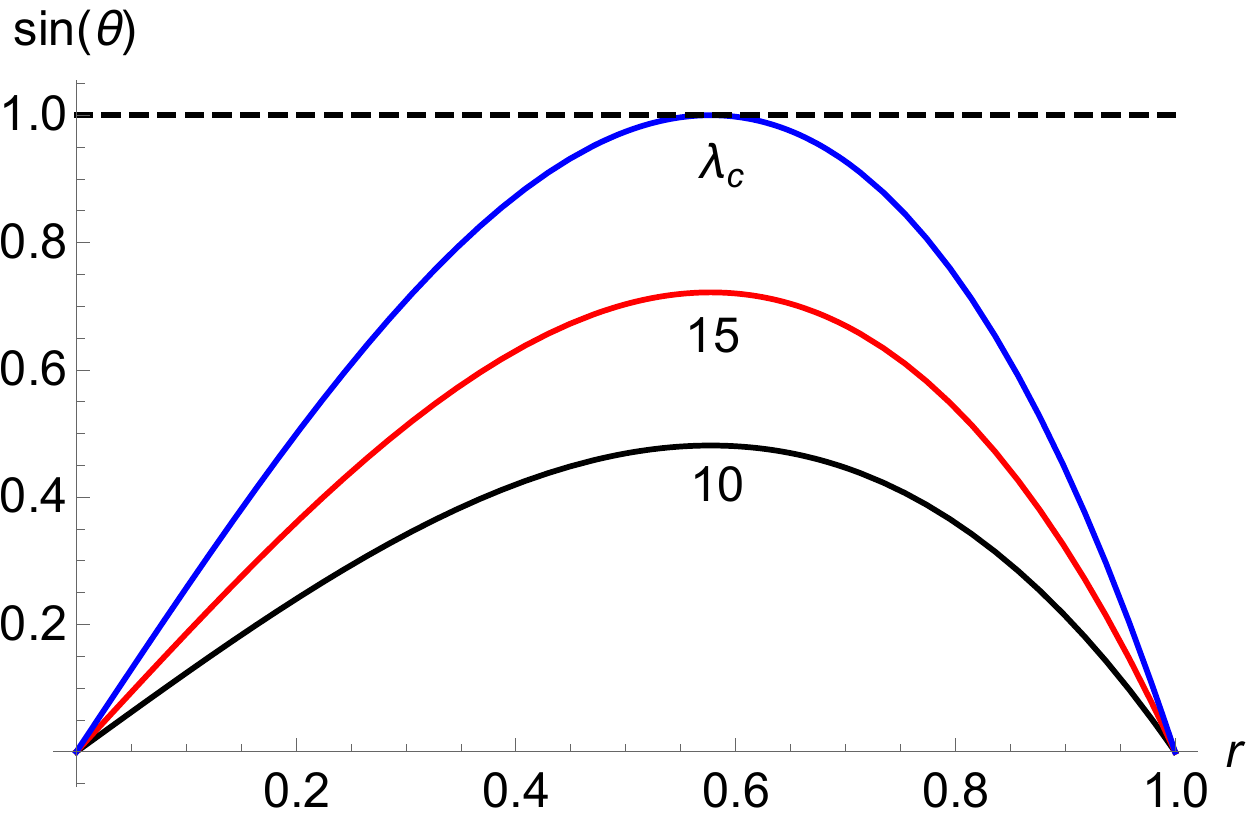}
\caption{Sine of the slope angle, $\theta$ (defined in fig. \ref{shapesschem}), vs. $r$ for $\lambda=10, 15$ and $\lambda_c=12\sqrt{3}$. The contact angle $\alpha=\pi/2$, i.e., $\theta(1)=0$. When $\lambda<\lambda_c$, the maximum slope corresponds to an inflection point with slope angle $0<\theta<\pi/2$.}
\label{sintheta}
\end{center}
\end{figure}
Fig. \ref{sintheta} shows $\sin\theta$ vs. $r$ as given by eq. (\ref{sint}) for various $\lambda$ and contact angle $\alpha=\pi/2$. Qualitatively similar shapes and trends arise for shapes with arbitrary $\alpha$. When $\lambda<\lambda_c$, $\theta$ attains a local maximum smaller than $\pi/2$; hence, the maximum slope angle coincides with an inflection point. However, at $\lambda_c$ the local maximum is $\theta_{\text{max}}=\pi/2$, and the axial meniscus length diverges with infinite slope and with a straight cylindrical shape of radius equal to the location where $\theta_{\text{max}}=\pi/2$; this is in agreement with previous studies, \cite{seebold1965,preziosi1987}. To find this location, we require 
\beq
\sin\theta=1,~~~~~\frac{d\sin\theta}{dr}=0,
\eeq
and obtain
\beq
\lambda_c=\frac{4}{r_c^3},~~~~
r_c=\frac{1}{2\cos\left(\frac{1}{3}(\pi-\alpha)\right)}.
\label{criticalrlambda}
\eeq
For $\alpha=\pi/2$, eqs. (\ref{criticalrlambda}) yield $\lambda_c=12\sqrt{3}$ and $r_c=1/\sqrt{3}$ as in fig \ref{sintheta}. Eqs. (\ref{criticalrlambda}) provide exact relations in support of the numerical results of Seebold \cite{seebold1965} for zero gravity. 
In particular, if one views the vertical slope location, $r_c$, as the radius of a straight circular cylinder, the rotational Bond number can be written as  $\overline\lambda_c\equiv r_c^3\lambda_c= 4$, providing a result that agrees with previous work \cite{seebold1965,preziosi1987}.\\
\\
Having shown that eq. (\ref{sint}) may be used both for a meniscus spanning the container with arbitrary contact angle $\alpha$ and for a wholly immersed bubble by setting $\alpha=0$, in the rest of this section we examine the case of $\alpha=\pi/2$, i.e., normal contact at the container wall. No generality is lost by focusing on this special case; on the contrary, it can be shown that, far from being special, critical shapes for arbitrary $\alpha$ may be obtained from the shape at any other $\alpha$ by suitable scaling manipulations.
We now show how $\alpha=\pi/2$ generates a master shape from which all other critical shapes with $0\le\alpha\le\pi$ may be constructed. The method is as follows: Starting with the critical shape for $\alpha=\pi/2$, we set $\lambda=\lambda_c=12\sqrt{3}$ and, using eq. (\ref{sint}), determine $r_w$ such that $\theta=\pi/2-\alpha$, i.e., $r_w$ is the location where the critical master shape has the same $\theta$ as the wall contact slope angle of interest. Setting $\sin(\pi/2-\alpha)=\cos\alpha=(\lambda_c/8) (r_w-r^{3}_w)$, we find:
\beq
r_w=\frac{2}{\sqrt{3}}\cos\left(\frac{1}{3}(\pi-\alpha)\right).
\label{rweq}
\eeq

We then rescale $r$,
\beq
r^* =\frac{r}{r_w}.
\eeq
Substituting for $r$ in eq. (\ref{rweq}) and manipulating to expose the binomial $(r^*-r^{*^3})$, we obtain the critical shape with $\theta=\pi/2-\alpha$ at $r=1$:
\beq
\sin\theta=r^*\cos\alpha+\frac{\lambda_c^*}{8} (r^*-r^{*^3}),
\label{sinthetacrit}
\eeq
where $\lambda_c^*$ satisfies eqs. (\ref{criticalrlambda}). This is the same equation that is solved in eq. (\ref{sint}) with arbitrary $\alpha$ for $\lambda=\lambda_c$ and demonstrates the generality of the $\alpha=\pi/2$ result for critical shapes. Even though full meniscus shapes do not exist at $\lambda=\lambda_c$ (since the axial meniscus length $H$ is infinite there), it is possible to use equation (\ref{sinthetacrit}) to predict portions of the meniscus shape, for values of $r$ not equal to $r_c$.

\subsection{Meniscus spanning the container radius}\label{infinitecontainer}
When the meniscus spans the cylinder radius $R$, we identify $d=R$ in eq. (\ref{lambda}). One integration of eq. (\ref{sint}) with a volume constraint fixes the absolute height of the interface. In a reference frame where the liquid volumes above and below $z=0$ are equal, the volume condition implies:
\beq
\int_0^1 h\,r dr = 0.
\label{volbc}
\eeq
The computation of $h(r)$ is performed numerically by integrating $h'(r)$ given by eq. (\ref{sint}). We focus on the case of normal contact ($\alpha=\pi/2$). This case is special only because the maximum slope location is $r_c=1/\sqrt{3}$ for all $\lambda$. Apart from this distinction, interface shapes are qualitatively similar when contact is not normal; and, as stated in Sec. \ref{properties}, interface shapes at criticality are easily scaled across different contact angles.  
\begin{figure}[h]
\begin{center}
\includegraphics[width=4in]{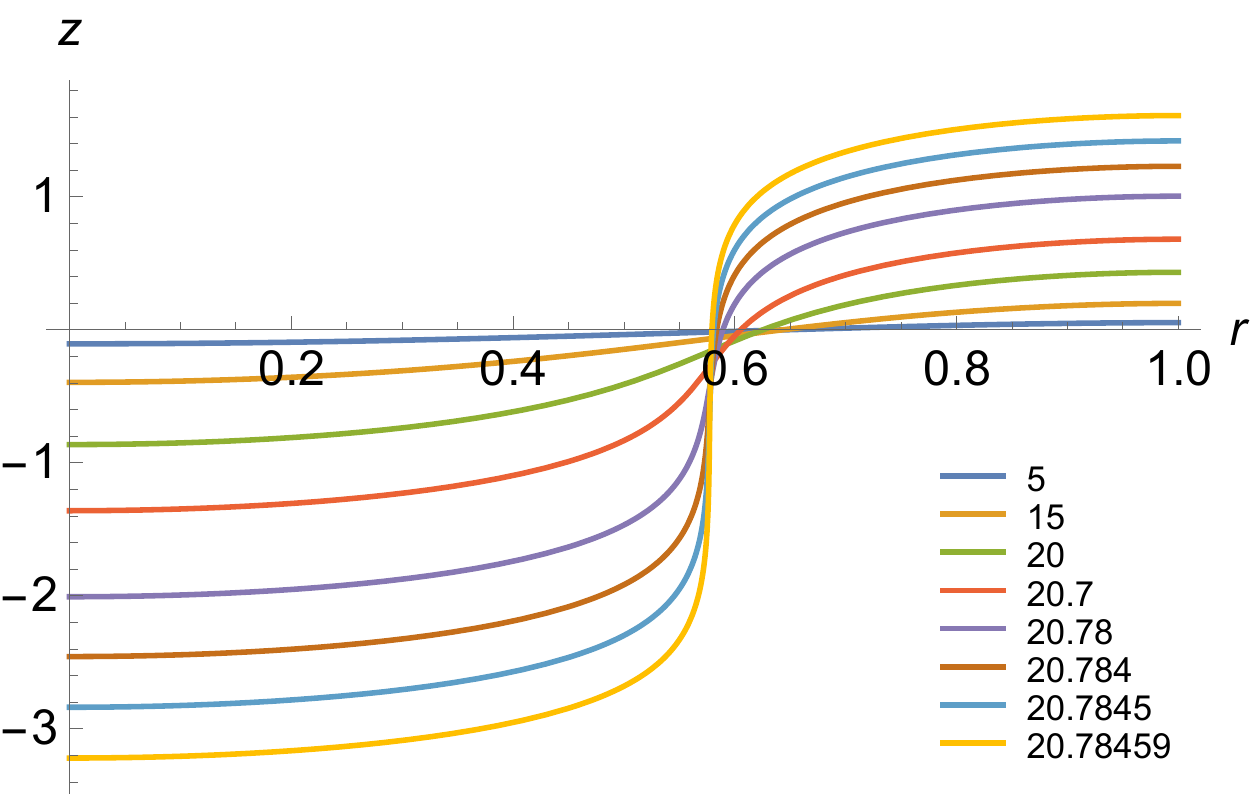}
\caption{Shapes with normal contact at $r=1$ (contact angle $\alpha=\pi/2$) for values of $\lambda$ noted in legend. The shape nearest $z=0$ is for $\lambda=5$. Consecutive values of $\lambda$ apply to shapes further away from $z=0$, showing the divergence of the meniscus axial length as $\lambda\to\lambda_c=12\sqrt{3}\approx 20.78459...$.}
\label{shapes90deg}
\end{center}
\end{figure}

\subsubsection{Meniscus asymptotics for $\lambda\to\lambda_c=12\sqrt{3}$}\label{infiniteasympt}
Fig. \ref{shapes90deg} shows interface shapes $z=h(r)$ when $\alpha=\pi/2$. It is clear that the axial meniscus length, $H$, defined in eq. (\ref{Hdef}), diverges as $\lambda\to\lambda_c$. Identifying the leading asymptotic behavior of this divergence is of considerable theoretical and practical interest because it is relevant to control devices such as rotating reactors where two immiscible fluids of differing densities are present. Since $h'(r_c)\to\infty$ as $\lambda\to\lambda_c$, it follows that $h'$ develops a narrowing peak around the maximum slope location, $r_c$; and that the area under the peak --though divergent-- depends to leading order on the shape of this peak only, i.e., it is independent of the details away from the peak. To begin, we identify the radial scale around the peak. Let
\begin{align}
\epsilon\equiv\lambda_c-\lambda,~~~~\eta\equiv\frac{r-r_c}{\epsilon^p},
\label{rlambda}
\end{align}
where $\epsilon\ll 1$, $\eta$ is a stretched radial distance centered at the peak, and $p$ needs to be determined. Approximating $h'$ from eq. (\ref{sint}) as $h'\sim \frac{1}{\sqrt{1-\sin^2\theta}}$ near $r_c$, and substituting $r$ and $\lambda$ from eq. (\ref{rlambda}), we find that
\beq
1-\sin^2\theta \sim \epsilon \frac{\sqrt{3}}{18}+\epsilon^{2p}\,9\,\eta^2,~~~\epsilon\to 0.
\eeq
This suggests that, for $h'$ to be integrable at $\eta=0$, we must have $p=1/2$ so that, to leading order as $\epsilon\to0$,
\beq
h'\sim h'_{\text{asy}}= \frac{1}{3\,\sqrt{\epsilon}\,\sqrt{\frac{1}{54\sqrt{3}}+\eta^2}}.
\label{hprimeasy}
\eeq
Anticipating the presence of an $O(1)$ constant following the leading divergent behavior as $\epsilon\to0$, we write:
\beq
H=\int_{-\frac{1}{\sqrt{3\epsilon}}}^{(1-\frac{1}{\sqrt{3}})\frac{1}{\sqrt{\epsilon}}}
\frac{d\eta}{3\,\sqrt{\frac{1}{54\sqrt{3}}+\eta^2}}+
\int_0^1 (h'-h'_{\text{asy}}) dr,
\eeq
which is equivalent to the first eq. (\ref{Hdef}) and where in the second integral $\eta$ has been written in terms of $r$ in the expression for $h'_{\text{asy}}$. 
The first integral may be evaluated in closed form and the second one is evaluated numerically. In the limit as $\epsilon\to0$, we obtain:
\begin{figure}[t]
\begin{center}
\includegraphics[width=4in]{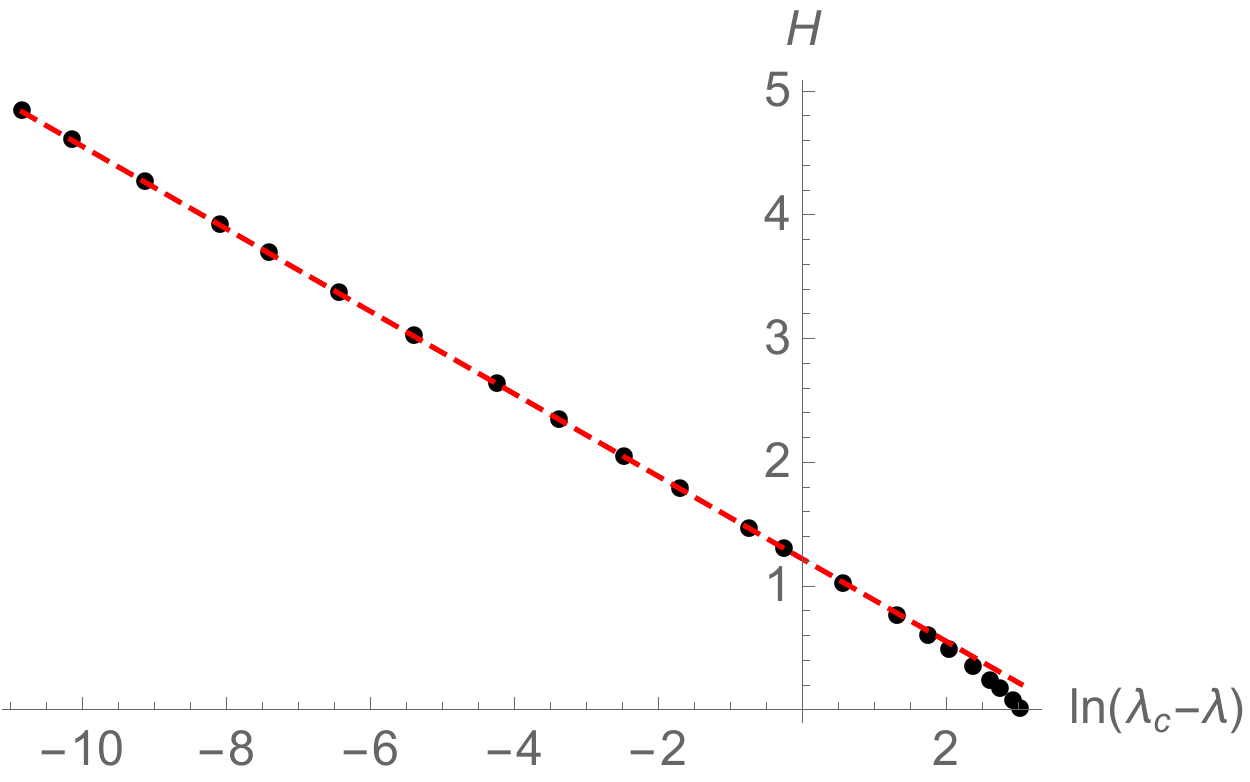}
\caption{Comparison of numerical and asymptotic evaluations of axial meniscus length $H$.  $\lambda_c-\lambda=\epsilon$ and $\lambda_c=12\sqrt{3}$. Black dots: numerical integration, eq. (\ref{Hdef}), with $h'$ from eq. (\ref{sint}). Dashed line: eq. (\ref{Hasyeq}). }
\label{Hasy}
\end{center}
\end{figure}
\begin{subequations}
\begin{align}
\int_{-\frac{1}{\sqrt{3\epsilon}}}^{(1-\frac{1}{\sqrt{3}})\frac{1}{\sqrt{\epsilon}}}
\frac{d\eta}{3\,\sqrt{\frac{1}{54\sqrt{3}}+\eta^2}} &\sim -\frac{1}{3}\ln\epsilon+\frac{1}{3}\ln[72(3-\sqrt{3})]+o(1)\\
\int_0^1 (h'-h'_{\text{asy}}) &\sim -0.2864+o(1),~~~~~~~\epsilon\to0
\end{align}
\label{H0eval}
\end{subequations}
Combining these results,
\beq
H\sim H_{\text{asy}}=-\frac{1}{3}\ln\epsilon+H_0,
\label{Hasyeq}
\eeq
where $H_0\approx1.218$. Fig. \ref{Hasy} shows the agreement between the numerical and asymptotic evaluations of the meniscus axial length.

\subsubsection{Series solution for $H(\lambda)$}\label{seriessolmeniscus}
In order to describe $H(\lambda)$ over the rest of the $\lambda$-domain, we now seek a series solution for $H(\lambda)$ in powers of $\lambda$ about $\lambda=0$. Since $h'(r,\lambda)=\tan(\theta)$, it follows that
\beq
h'(r,\lambda)=\frac{\sin\theta}{\sqrt{1-\sin^2\theta}}.
\label{hprime}
\eeq
For $\alpha=\pi/2$, this expression may be expanded as a Taylor series about $\lambda=0$ using eq. (\ref{sint}) as:
\beq
h'(r,\lambda)=\sum_{n=0}^{\infty} a_{2n+1}(\sin\theta)^{2n+1}=
\sum_{n=0}^{\infty} a_{2n+1}\left(\frac{r-r^3}{8}\right)^{2n+1} \lambda^{2n+1},
\label{hprimeseries}
\eeq
where
\begin{subequations}
\begin{align}
a_1=1,~~~&a_{n}=a_{n-2}\left(\frac{n-2}{n-1}\right),~~~n=3,5,7,9...
\label{an1}\\
&a_n=0,~~~~~n=0,2,4,...
\label{an2}
\end{align}
\end{subequations}
As seen in fig. \ref{sintheta}, $\sin\theta$ remains below 1 for $\lambda<\lambda_c$. Thus, the series in eq. (\ref{hprimeseries}) converges for $\lambda<\lambda_c$ and may be integrated term by term to obtain the meniscus length:
\beq
H(\lambda)=\sum_{n=0}^{\infty} a_{2n+1}\, b_{2n+1} \lambda^{2n+1},
\label{HTaylor}
\eeq
where
\begin{subequations}
\beq
b_n=\frac{1}{8^n}\int_0^1(r-r^3)^n dr,~~~~~n~~\text{odd},
\eeq
\beq
b_n=0,~~~~~n~~\text{even}.
\eeq
\end{subequations}
For $n$ odd, it is easy to show that $b_1=1/32$ and
\beq
b_{n+2}=\frac{1}{64}
        \frac{(n+2)(n+1)\left(\frac{n+1}{2}\right)}
             {\left(\frac{3n+7}{2}\right)\left(\frac{3n+5}{2}\right)
              \left(\frac{3n+3}{2}\right)}~b_n.
              \label{bn2}
\eeq
The ability to compute all the terms of the infinite series permits evaluation of the radius of convergence of the series given in eq. (\ref{HTaylor}). The ratio test guarantees convergence iff
\beq
\lim_{n\to\infty}\frac{a_{2n+1}}{a_{2n-1}}\frac{b_{2n+1}}{b_{2n-1}}\lambda^2<1.
\eeq
Evaluation of this criterion using eqs. (\ref{an1}) and (\ref{bn2}) shows that the series does converge for $\lambda<12\sqrt{3}=\lambda_c$ as was stated above. Thus, eq. (\ref{HTaylor}) is an exact solution. Unfortunately, though, the convergence is poor and nonuniform with increasing $\lambda$ beyond $\lambda\approx 15$ due to the influence of a logarthmic singularity at $\lambda=\lambda_c$, see Sec. \ref{infiniteasympt}. In Sec. \ref{approximants} we use asymptotic approximants to generate a rapidly converging and uniform representation of $H(\lambda)$, defined in eq. (\ref{Hdef}), over the entire range $0\le\lambda< \lambda_c$.

\subsection{Spinning bubble}
If the gas volume in a finite container is small enough, or the axial dimension of the container is long enough, the bubble can become arbitrarily close to critical (i.e., it can adopt a nearly straight circular cylindrical shape with locally curved ends) before the interface touches the end plates of the container. This is the basis for the well-known spinning bubble method to determine surface tension \cite{vonnegut1942, joseph1994}. In this geometry (see fig. \ref{bubbleschem}) and in zero-gravity, bubbles exist for ${\lambda<\lambda_c}$$=4$, and the characteristic length, $d$, is the maximum bubble radius, $R_b$. For extensive detail on the challenges of interpreting and operating the spinning drop tensiometer in a gravitational field, see Manning \& Scriven \cite{scriven1977} and references therein.
\\
\\
Using the bubble maximum radius, $R_b$, as the characteristic length, the slope angle of the interface relative to the $r$-axis is found by integrating eq. (\ref{sint}) subject to $\theta(0)=0$, $\theta(1)=\pi/2$:
\beq
\sin\theta=r+\frac{\lambda}{8}r(1-r^2).
\label{spinshape}
\eeq

\begin{figure}[t!]
\begin{center}
\includegraphics[width=3.5in]{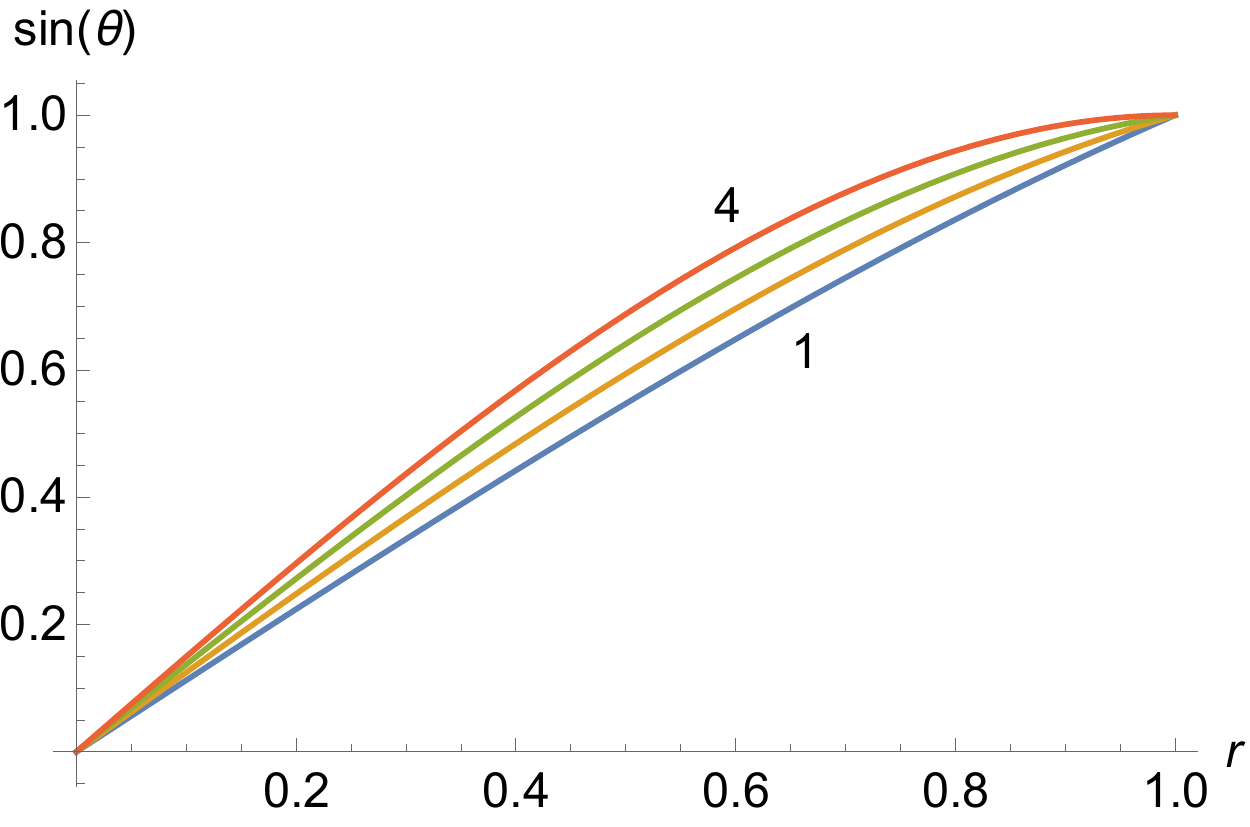}
\caption{Shapes of spinning bubbles at various $\lambda=1,2,3,4$. Only $\lambda=4$ has $d(\sin\theta)/dr=0$ at $r=1$ .}
\label{sintehetabubble}
\end{center}
\end{figure}
Since $\sin\theta=h'/\sqrt{1+h'^2}$, numerical integration of $h'$ subject to $h(0)=0$ yields the interface shape, $z=h(r)$. The bubble has infinite slope ($\theta_{\text{max}}=\pi/2$) for all $\lambda$ at $r=1$. This contrasts the meniscus analyzed in Sec. \ref{infinitecontainer} where $\theta_{\text{max}}<\pi/2$ and coincides with an inflection point located in $0<r<1$. In this case, $\theta_{\text{max}}\to\pi/2$ \textit{only} as $\lambda\to\lambda_c$. To probe the character of the spinning bubble shape, we note that, when $\lambda<4$, a point located on the bubble's equator, denoted by $(r,\bar{z})=(1,0)$, corresponds to a local maximum of $r$. Therefore, $dr/d\bar{z}=0$ and $d^2r/d\bar{z}^2<0$ there. From geometrical considerations, this implies that $1-r\sim A \bar{z}^2$, so that $\sin\theta\sim 1-2A(1-r)$ as $r\to1$ for some constant $A>0$. In contrast, when the shape is critical at $\lambda=\lambda_c$, the end-cap shape approaches a straight cylinder asymptotically at a distance from the bubble tip that is large compared with the radius; therefore, in the critical condition, $1-r\sim\exp(-Bz)$ for some constant $B>0$, where $(r,z)=(0,0)$ is the tip location, $z\gg1$ and the bubble is in $z>0$. This implies that $\sin\theta\sim 1-(B^2/2)(1-r)^2$ as $r\to1$ for $\lambda=\lambda_c$. We conclude, therefore, that the critical shape requires
\beq
\frac{d(\sin\theta)}{dr}=0~~~~\text{at}~~~r=1,~~\lambda=\lambda_c=4,
\eeq
whereas subcritical shapes satisfy
\beq
\frac{d(\sin\theta)}{dr}>0~~~~\text{at}~~~r=1,~~\lambda<4.
\eeq
Not surprisingly, the shapes of eq. (\ref{spinshape}), a few of which are shown in fig. \ref{sintehetabubble}, display these properties. In contrast to the meniscus spanning the cylinder radius, the distinct character of the spinning bubble configuration is that $\sin\theta=1$ always at $r=1$; but $d(\sin\theta)/dr\ne0$ at $r=1$ unless $\lambda=\lambda_c=4$.

In the rest of this section we examine the bubble axial length, $H(\lambda)$. We derive the asymptotic behavior as $\lambda\to\lambda_c$ and develop an exact solution as a series in powers of $\lambda$. Both analyses can be used to better inform the quality of the critical character of a bubble formed in an experiment.

\begin{figure}[t!]
\begin{center}
\includegraphics[width=3.5in]{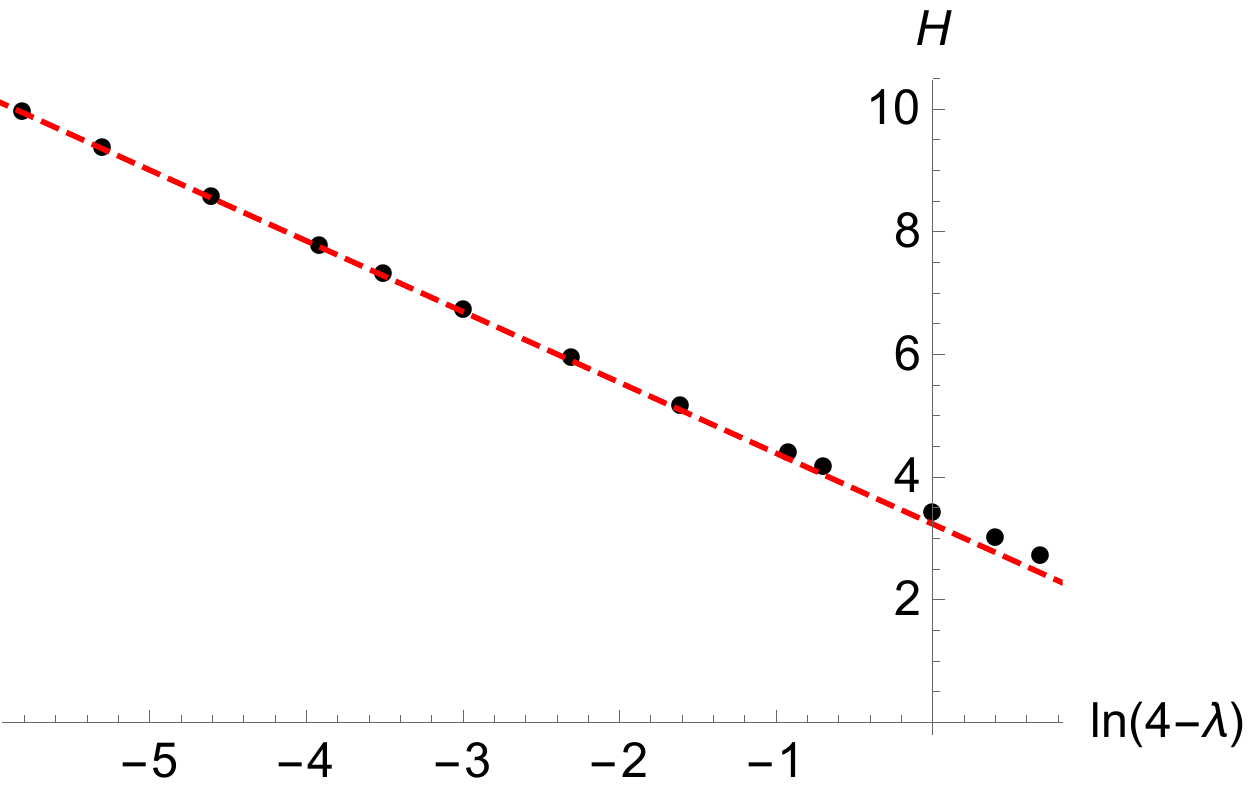}
\caption{Comparison of numerical and asymptotic evaluations of axial meniscus length, $H$, for the spinning bubble vs. $\lambda$. Black dots: numerical integration, eq. (\ref{Hdef}). Dashed line: Asymptotics, eq. (\ref{Hasyspinning}).}
\label{criticalbubble}
\end{center}
\end{figure}
\subsubsection{Asymptotics for $\lambda\to\lambda_c= 4$}
We derive the asymptotic behavior of the shape as $\lambda\to\lambda_c=4$ by the same method of Sec. \ref{infiniteasympt}. Let $\epsilon\equiv4-\lambda$ and $\eta\equiv(r-1)/\epsilon^p$, where $p>0$ is to be determined. It may be shown from eq. \ref{spinshape} that, when $\epsilon$ is small,
\beq
\sin^2\theta\sim 1-3\epsilon^{2p}\eta^2+\epsilon^{1+p}\frac{\eta}{2}+...
\eeq
Approximating $h'$ near the peak as $1/\sqrt{1-\sin\theta^2}$, the only choice that ensures integrability of $h'$ at $\eta=0$ is $p=1$, yielding
\beq
h'\sim h'_{asy}=\frac{1}{\epsilon}\frac{1}{\sqrt{3\eta^2-\frac{\eta}{2}}}.
\eeq
Using the same methodology as in Sec. \ref{infiniteasympt}, we obtain
\beq
H_{asy}=\sim -\frac{2}{\sqrt{3}}\ln(4-\lambda)+3.2332, ~~~~\lambda\to4.
\label{Hasyspinning}
\eeq
The (dimensionless) bubble volume depends on $\lambda$:
\beq
V(\lambda)\equiv\frac{\widetilde V}{R_b^3}=2\left[\pi h(1,\lambda)
-2\pi\int_0^1 h(r,\lambda)\,r\,dr\right]
\label{dimlessvol}
\eeq
and provides a relation between the dimensional volume, $\widetilde V$, and maximum radius $R_b$. The closer the bubble is to the critical configuration the closer its shape is to a straight cylinder of dimensionless radius 1; hence, the dimensionless volume grows progressively more linearly with $\pi H$ as $\lambda\to\lambda_c=4$. The limit of eq. (\ref{dimlessvol}) as $\lambda\to4$ is
\beq
V\sim \pi H - 4.18,~~~~H\gg1,
\label{vasy}
\eeq
which is in good agreement with the asymptotic behavior of Ross' exact expression for the volume (eqn. 15 in ref. \cite{ross1968}) as $\lambda\to \lambda_c=4$; refer to figs. 7 and 8. 

\begin{figure}[t!]
\begin{center}
\includegraphics[width=3.5in]{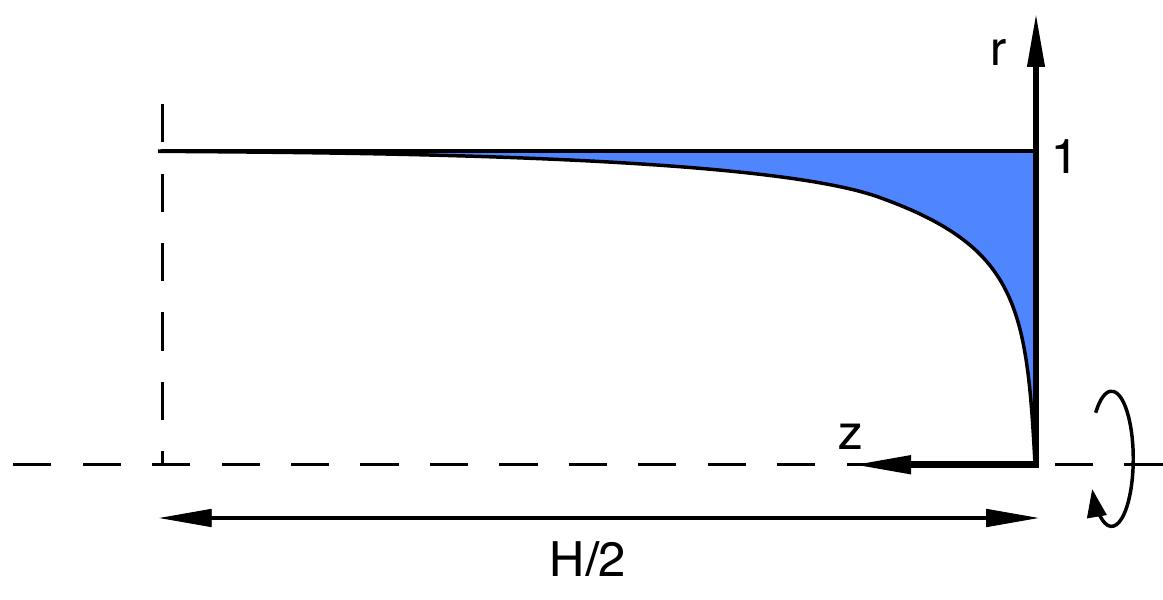}
\caption{Schematic of a half-bubble spinning about the $z-$axis. The blue area represents the liquid and is equal to the volume deficit of the actual bubble relative to a straight circular cylinder where the bubble is inscribed (see eq. \ref{vasy}). The spinning container is larger than the bubble size and is not shown.}
\label{spinningbubble}
\end{center}
\end{figure}

\begin{figure}[t!]
\begin{center}
\includegraphics[width=3.5in]{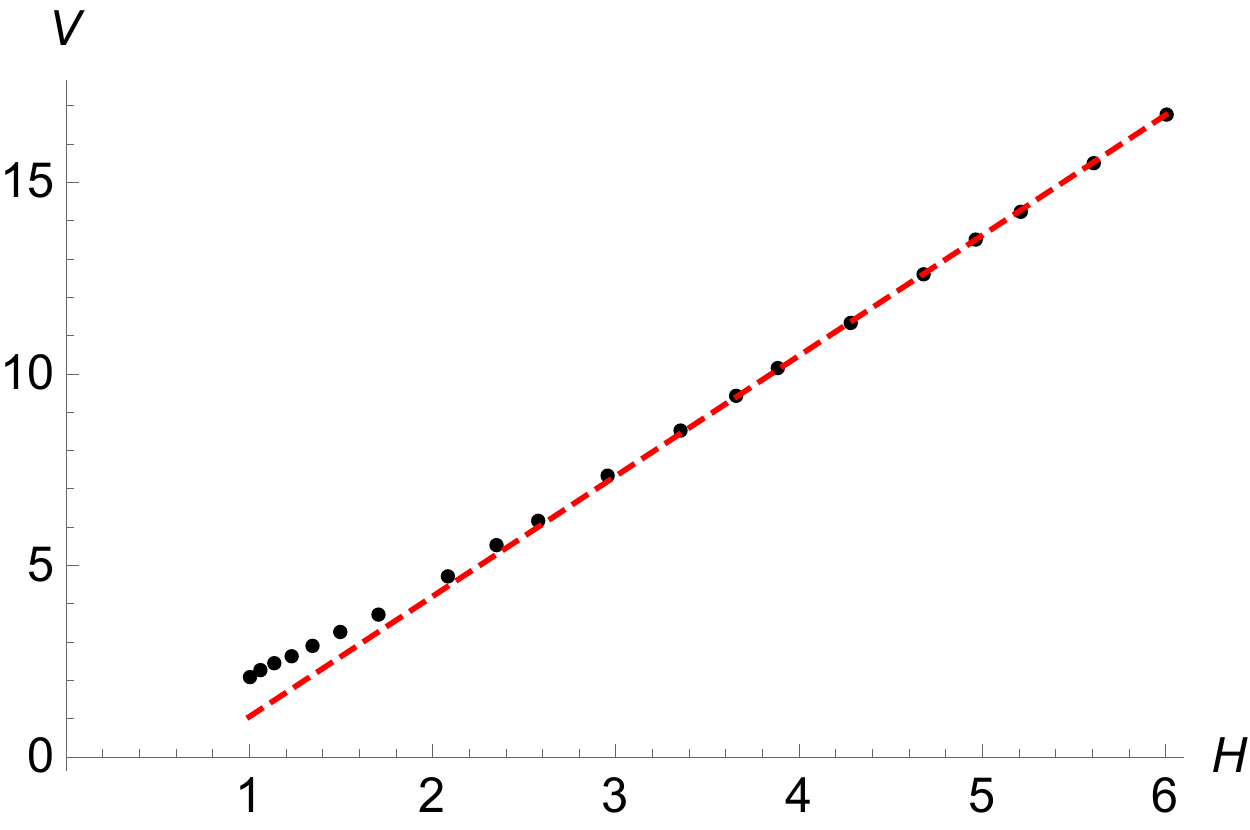}
\caption{Dimensionless volume, $V$, vs bubble length $H$. Dashed line: equation (\ref{vasy}).}
\label{VvsH}
\end{center}
\end{figure}

\subsubsection{Series solution for $H(\lambda)$}\label{seriesbubble}
In order to construct a Taylor series representation of the bubble axial length, defined in eq. (\ref{Hdef}), we note that, as in Sec. \ref{infinitecontainer}, $h' = \sin\theta/\sqrt{1-\sin^2\theta}$. However, in contrast to that analysis, here $\sin\theta\ne O(\lambda)$, which complicates evaluation of a Taylor series in powers of $\lambda$ for $H(\lambda)$. Let us first generate the series
\beq
h'(r,\lambda)=\sum_{n=0}^{\infty} c_n(r) \lambda^n.
\label{seriesbubble1}
\eeq
Starting from eq. (\ref{sint}), we write
\beq
\sin\theta=A(r) + \lambda B(r),
\eeq
where $A=r$ and $B=\frac{r-r^3}{8}$. It follows that the denominator in the expression for $h'(=\tan\theta)$ is
\beq
\sqrt{1-A^2-2\lambda AB-\lambda^2 B^2}=
   \left(\sum_{j=0}^2 a_j \lambda^j\right)^{1/2},
   \label{hprimedenom} 
\eeq
where $a_0=1-A^2$, $a_1=-2AB$ and $a_2=-B^2$. Using J.C.P. Miller's formula for the series expansion of a series raised to any power \cite{henrici1956}, we evaluate the series for the inverse of (\ref{hprimedenom}):
\beq
\left(\sum_{j=0}^2 a_j \lambda^j\right)^{-1/2}=\sum_{j=0}^{\infty}b_j\lambda^j,
\label{miller}
\eeq
to find the following recursion for the coefficients,
\begin{align}
c_0&=A\,b_0 \nonumber\\
c_{n>0}&=A\,b_n+B\,b_{n-1},
\end{align}
where
\begin{align}
b_0&=\frac{1}{\sqrt{1-A^2}},~~~~~~b_1=\frac{AB}{\sqrt{1-A^2}}\,;\nonumber\\
b_{n>1}&=
        -\frac{1}{n(1-A^2)}\left[\left(\frac{1}{2}-n\right) 2A\,B\,b_{n-1}
       +(1-n)\,B^2\,b_{n-2}\right]. \label{miller2}
\end{align}
Since $c_n$'s are linear combinations of $b_n$'s, convergence properties of the series in eq. (\ref{seriesbubble1}) can be determined from those of $\sum_n b_n\lambda^n$. Dividing through by $b_{n-1}$ in eq. (\ref{miller2}) we form two ratios of consecutive $b_n$. Assuming that this ratio has a limit  as $n\to\infty$, denoted $Q_{\infty}(r)$, solution of a quadratic equation yields $Q_{\infty}(r)=(r^2+r)/8$. Based on the ratio criterion, convergence is guaranteed with $0\le r\le 1$ iff Max$_{r}[Q_{\infty}(r)] \,\lambda <1$, i.e., $\lambda<4$. As in the problem of Sec. \ref{infinitecontainer}, the series in eq. (\ref{seriesbubble1}) converges in the entire range of $\lambda$ where shapes exist, i.e., $0\le\lambda<4$, and is therefore an exact solution. It can therefore be integrated term-by-term to produce another convergent exact solution for $H(\lambda)$, defined in eq. (\ref{Hdef}), e.g.,
\beq
H(\lambda)=2\sum_{n=0}^{\infty}\left(\int_0^1 c_n(r) dr\right) \lambda^n=
           \sum_{n=0}^{\infty} C_n \lambda^n.
\label{HTaylor2}
\eeq
(In the appendix we show an explicit evaluation of $C_n$.) Because convergence of the series in eq. (\ref{HTaylor2}) is poor as $\lambda$ increases beyond $\lambda\approx 3.5$, in Sec. \ref{approximants} we show how to implement the method of asymptotic approximants to obtain an analytical expression for $H(\lambda)$ that is uniform across the entire range $0\le\lambda\le\lambda_c=4$.

\subsection{Approximants}\label{approximants}
Asymptotic approximants provide uniformly convergent approximations to the axial lengths $H(\lambda)$, as given by eq. (\ref{HTaylor}) for a rotating meniscus and eq. (\ref{HTaylor2}) for a spinning bubble, over the entire respective intervals $0\le\lambda\le\lambda_c$. Interested readers may consult Barlow et al. \cite{barlow2017} and references therein for an extensive presentation of the method applied to a wide range of problems in mathematical physics.\\
\\
Briefly, asymptotic approximants go beyond the well-known Pad\'e approximants  in that they incorporate asymptotic behaviors that are often singular in ways other than just poles \cite{benderorszag}, thus dramatically improving the approximant's power to extend the region of convergence. Both power series for $H(\lambda)$ in the present work (eqs. (\ref{HTaylor}) and (\ref{HTaylor2})) have a logarithmic divergence at their respective $\lambda_c$. Because we have the power series expanded about $\lambda=0$, as well as the logarithmic divergence behavior as $\lambda$ approaches $\lambda_c$, an asymptotic approximant may be used to join these behaviors. In the two problems considered here, we propose the following approximant for $H(\lambda$) defined in eqs. (\ref{HTaylor}) and (\ref{HTaylor2}):
\beq
H_A(\lambda, N) = \sum_{n=0}^{N} A_n (\lambda_c-\lambda)^n + A_L +
               B_L \ln(\lambda_c-\lambda),
\label{approximantform}
\eeq
where $A_L$ and $B_L$ have been computed from the respective asymptotic analyses in eqs. (\ref{Hasyeq}) and (\ref{Hasyspinning}). The coefficients $A_n$ are determined from the condition that the $N$-term Taylor series of $H_A(\lambda, N)$ about $\lambda=0$ is equal to the $N$-term Taylor series of $H(\lambda)$ in eqs. (\ref{HTaylor}) and (\ref{HTaylor2}). The form in eq. (\ref{approximantform}) imposes the asymptotic logarithmic divergence as $\lambda\to\lambda_c$.
\begin{figure}[h!]
\begin{center}
\includegraphics[width=3.5in]{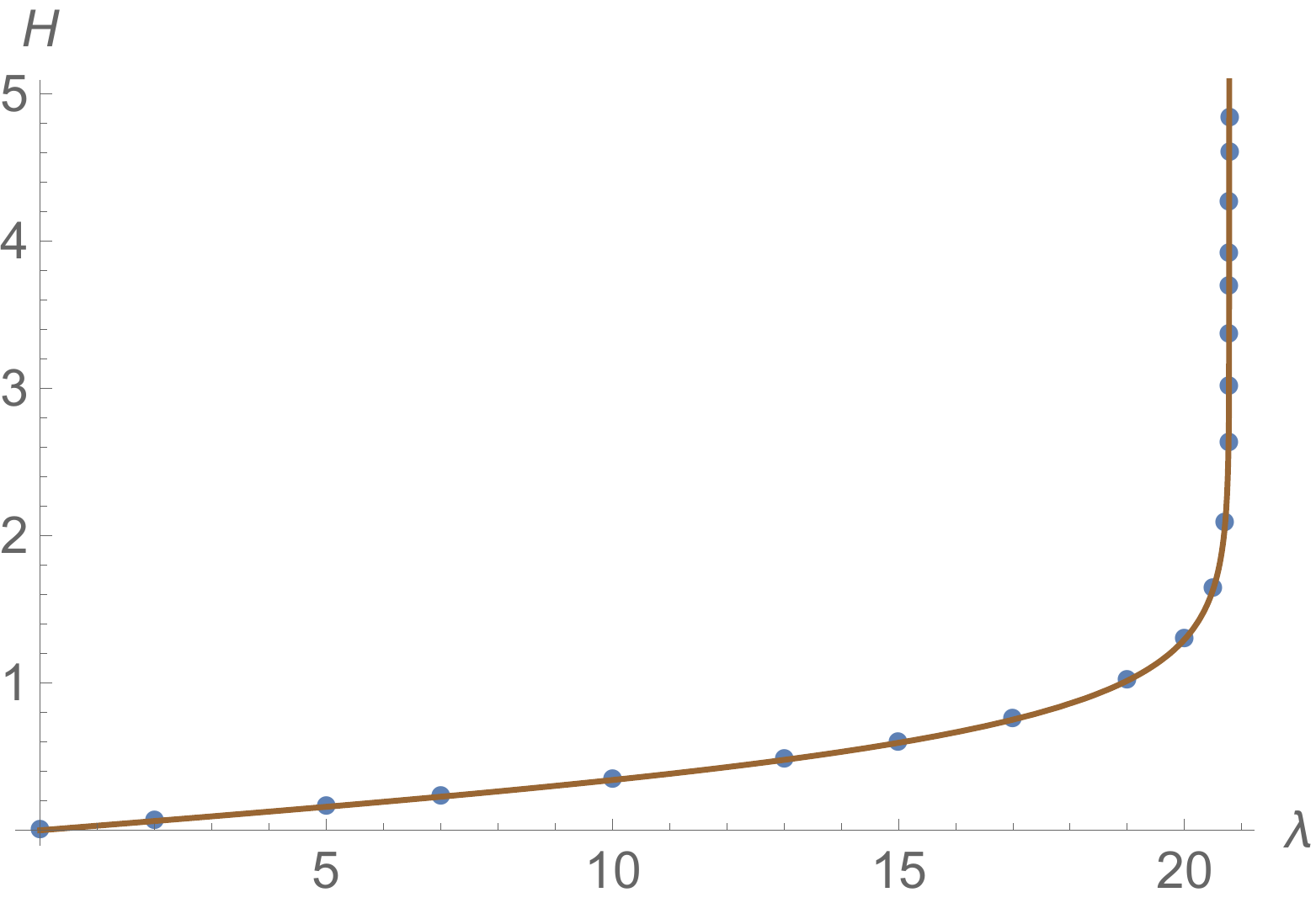}
\caption{The axial length of the rotating meniscus versus $\lambda$. Dots: numerical evaluation of $H$, eq. (\ref{Hdef}). Line: 20-term $H_A$ from eq. (\ref{approximantform}); $A_L=1.218, B_L=1/3$. At $\lambda=0$ the interface is flat, therefore $H(0)=0$.}
\label{HAMeniscus}
\end{center}
\end{figure}
\begin{figure}[h!]
\begin{center}
\includegraphics[width=3.5in]{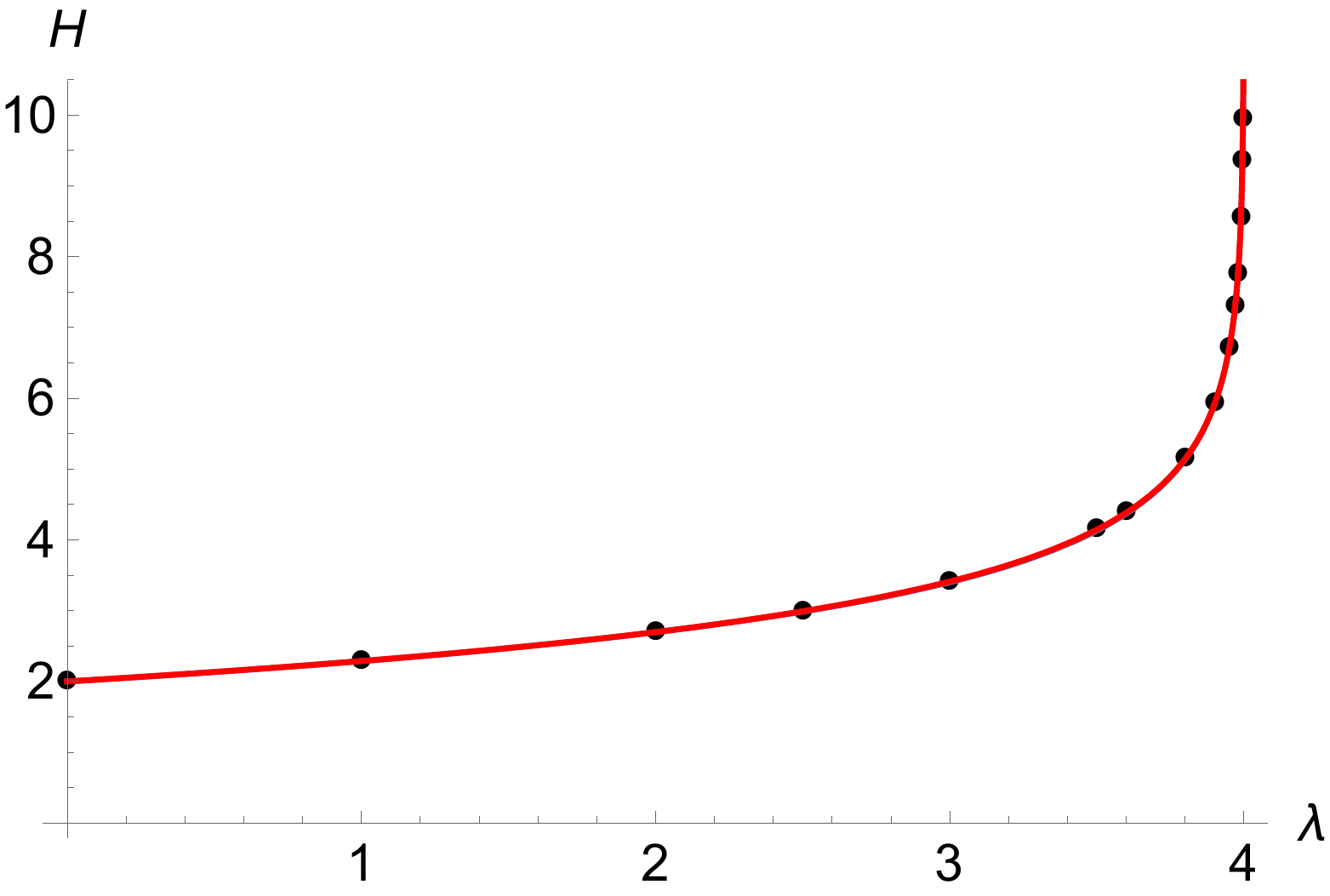}
\caption{The length of the spinning bubble versus $\lambda$. Dots: numerical valuation of $H$, eq. (\ref{Hdef}). Two solid lines, the 5-term and 10-term $H_A$, are indistinguishable; $A_L=3.2332, B_L=2/\sqrt{3}$. At $\lambda=0$ and in zero-gravity, the bubble is spherical, i.e., $H(0)=2$.}
\label{HABubble}
\end{center}
\end{figure}
\\
\\
The form of an asymptotic approximant to a given function is not uniquely determined but experience allows one to pose forms that exhibit superior convergence. In this work we have not attempted to optimize the form of the approximant that minimizes the number of terms $N$ required to produce a given error.  Figures \ref{HAMeniscus} and \ref{HABubble} show $H(\lambda)$ from the numerical integration of eq. (\ref{Hdef}) together with $H_A(\lambda, N)$, for the rotating meniscus and the spinning bubble, respectively, and for various number of terms, $N$, in the approximant.\\
\\
We define the error of the $N$-term approximant, $\Delta_{N}(\lambda)$, as the pointwise absolute error between $H_A(\lambda,N)$ and the numerical values of $H(\lambda)$. Both approximants seem to converge to the numerical calculation as $N$ increases; but, as shown in figs. \ref{errorvortex} and \ref{errorbubble}, the convergence has a small non-uniformity near $\lambda_c$. This is a well-known phenomenon due to remaining singularities that approach zero as $\lambda\to\lambda_c$ \cite{harkin2021_2}. The largest errors of the most accurate approximants calculated are $\Delta_{20}\approx4\times10^{-4}$ for the rotating meniscus and $\Delta_{15}\approx7.2\times10^{-3}$ for the spinning bubble. 
\begin{figure}[t!]
\begin{center}
\includegraphics[width=3.5in]{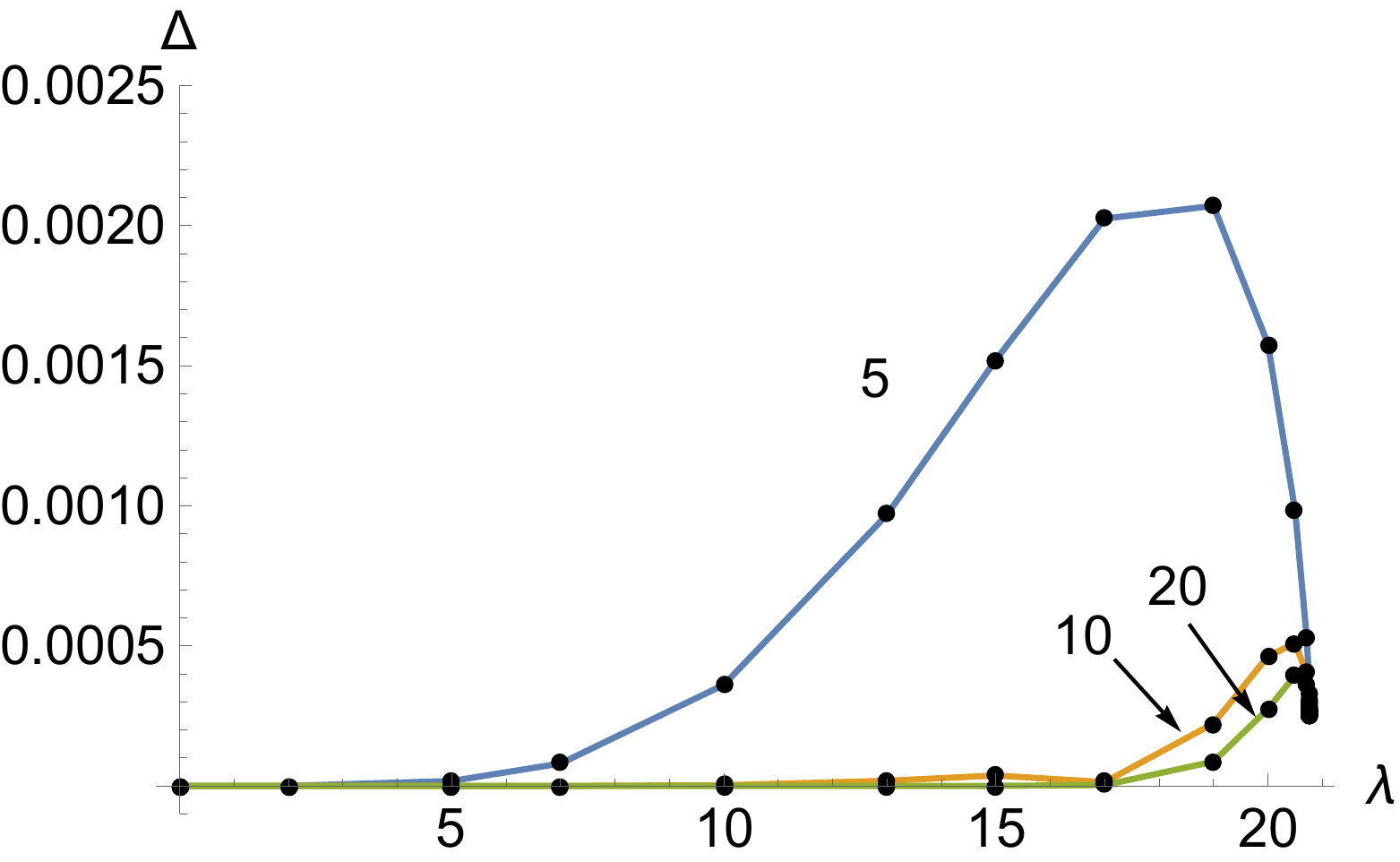}
\caption{Pointwise error, $\left|H(\lambda)-{H_A}(\lambda)\right|$, at numerically calculated points for the rotating meniscus for approximants with number of terms $N=5,10,20$. Lines are indicative only and are used as a guide to the eye.}
\label{errorvortex}
\end{center}
\end{figure}
\begin{figure}%[h!]
\begin{center}
\includegraphics[width=3.5in]{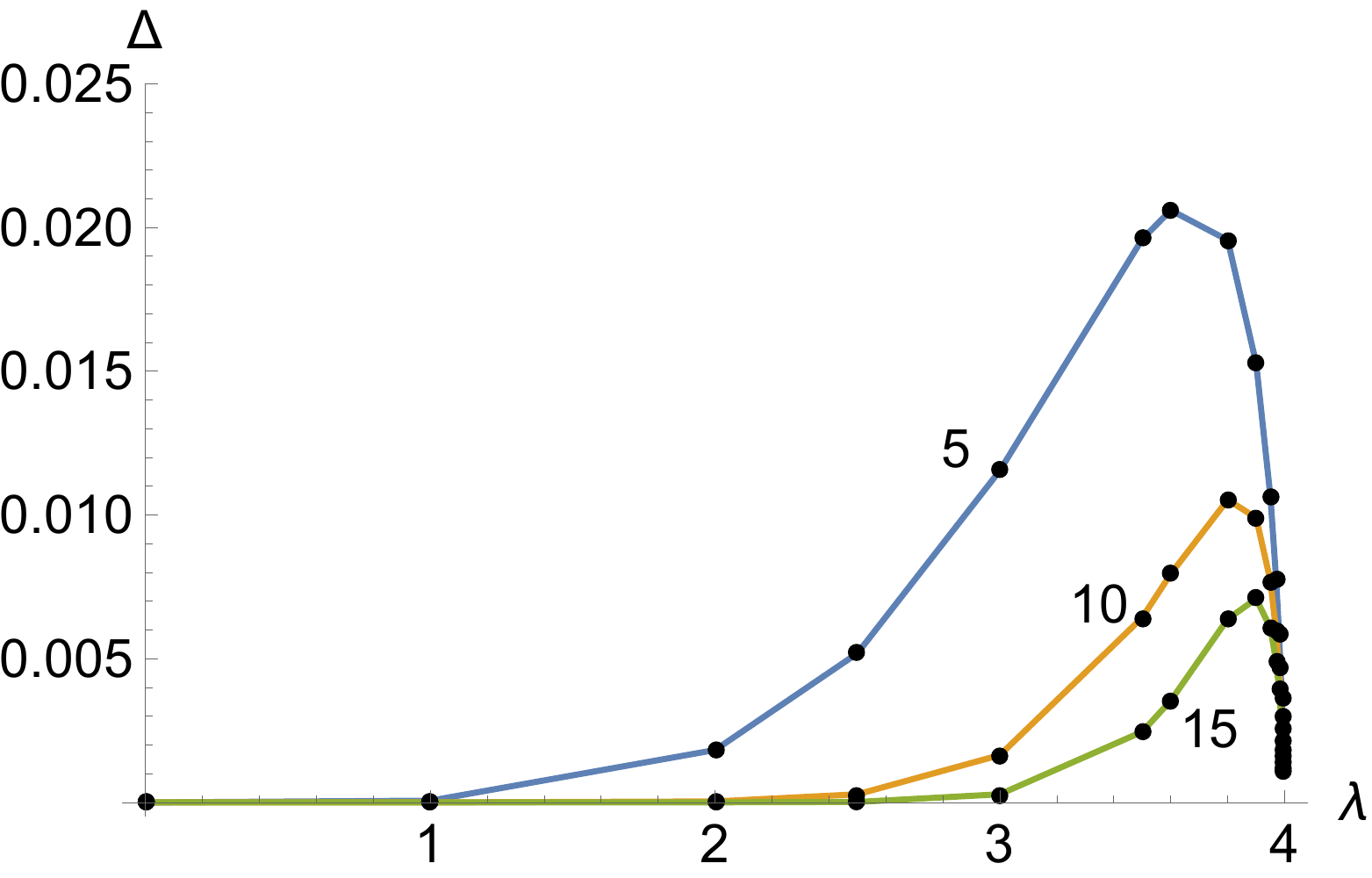}
\caption{Pointwise error, $\left|H(\lambda)-{H_A}(\lambda)\right|$, at numerically calculated points for the spinning bubble for approximants with number of terms $N= 5, 10, 15$. Lines are indicative only and are used as a guide to the eye.}
\label{errorbubble}
\end{center}
\end{figure}

\section{Discussion}\label{discussion}
\begin{figure}[t]
\begin{center}
\includegraphics[width=3.5in]{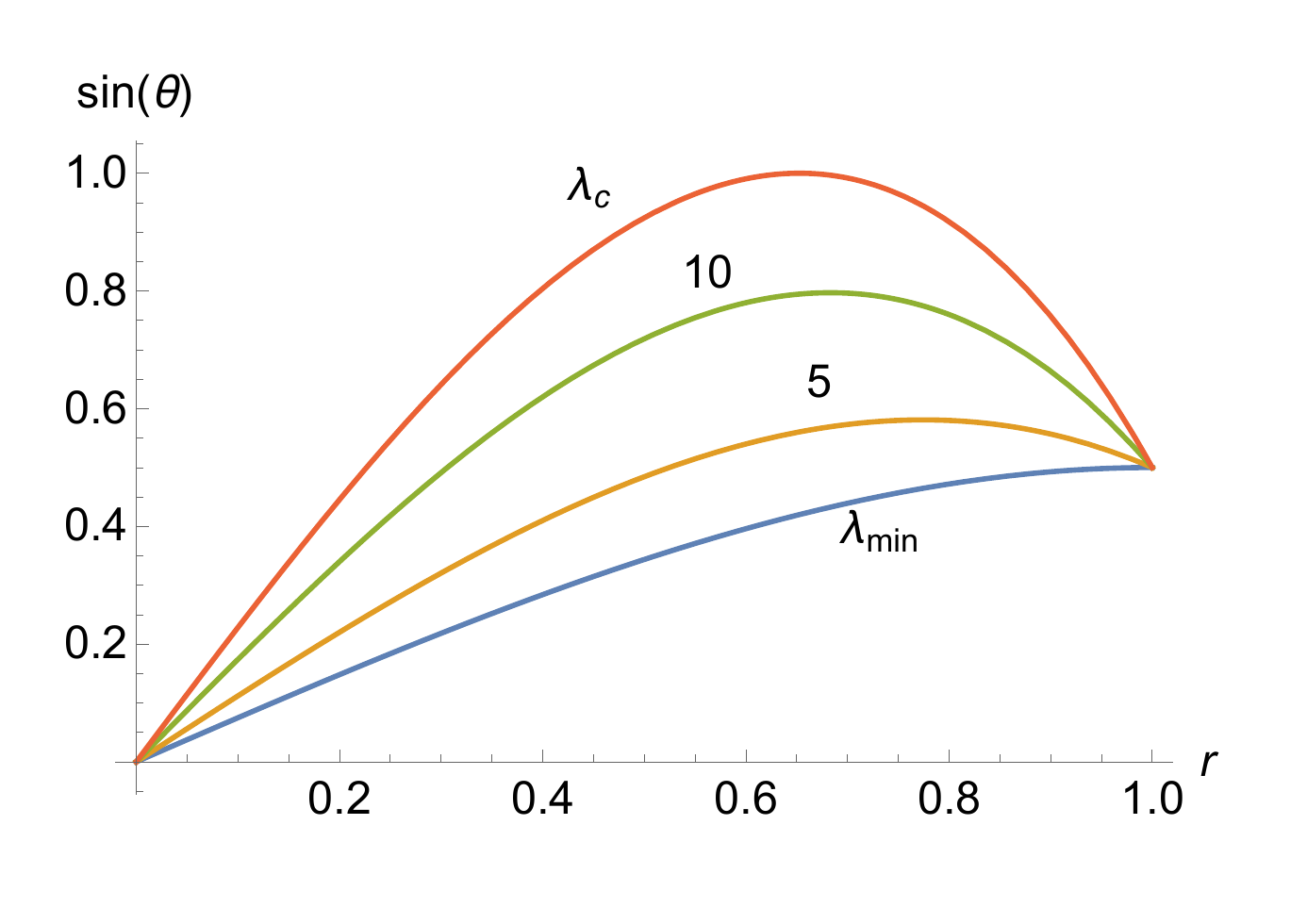}
\caption{Shapes for $\alpha=\pi/3$ (60-degree contact angle) for $\lambda=2, 5, 10$ and 14.39. $\lambda_{\text{min}}=2$ and $\lambda_c\approx 14.39$. The location of the maximum at each $\lambda$ is plotted below in fig. \ref{rmaxtheta}.}
\label{30degsint}
\end{center}
\end{figure}

\begin{figure}[t]
\begin{center}
\includegraphics[width=3.5in]{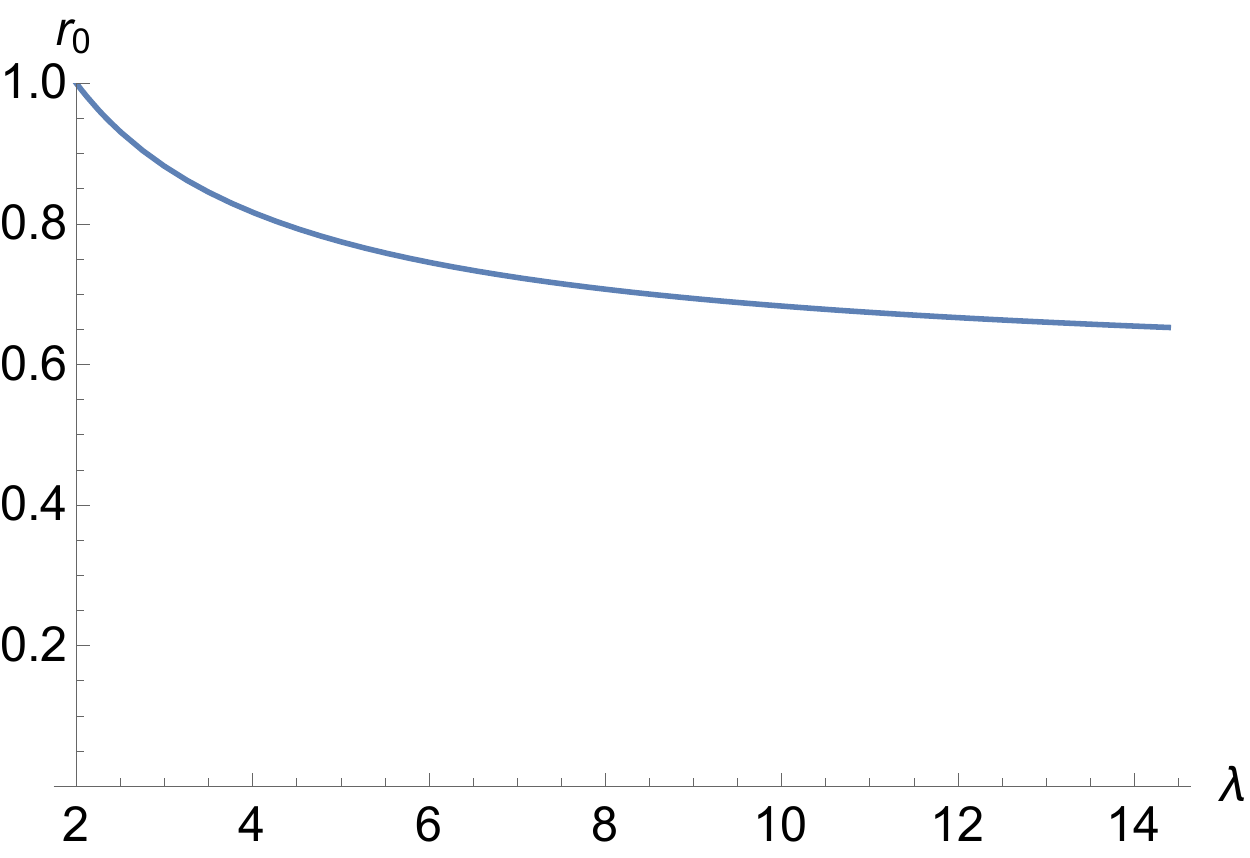}
\caption{$r_0$ vs. $\lambda$ between $\lambda_{\text{min}}=2$ and $\lambda_c\approx 14.385$, for $\alpha=\pi/3$ (60-degree contact angle). }
\label{rmaxtheta}
\end{center}
\end{figure}

In configurations with arbitrary wall contact angle, $\alpha$, the shape is described by eq. (\ref{sint}), rewritten here,
\beq
\sin\theta=r\cos\alpha +\frac{\lambda}{8}(r-r^3).
\eeq
Since, for the general cases of $\alpha\ne 0$, critical shapes arise from the progressive steepening of an oblique inflection point, it is instructive to find, for a given $\alpha$, $\lambda$ for which $\sin\theta$ has a maximum at $r=1$; this value of $\lambda$ is a lower bound for $\lambda_c$. By requiring that $d(\sin\theta)/dr=0$ at $r=1$ this value is found to be $\lambda_{\text{min}}=4\cos\alpha$. In order to achieve criticality in $0<r<1$, $\lambda$ must be greater than $\lambda_{\text{min}}$. The location of $\theta_{\text{max}}$ for $\lambda>\lambda_{\text{min}}$ is
\beq
r_0=\left[\frac{1}{3}\left(1+\frac{8}{\lambda}\cos\alpha\right)
    \right]^{1/2}.
\eeq
Inserting this result into the expression for $\sin\theta$ above, we find the maximum value of $\sin\theta$ for a given $\alpha$ as function of $\lambda$:
\beq
\sin\theta|_{r_0}=\frac{1}{3\sqrt{3}}\,\frac{\lambda}{8}
                  \left(1+\frac{8}{\lambda}\cos\alpha\right)^{3/2}.
\eeq
This expression becomes 1 when $\lambda=\lambda_c$ as given in eq. (\ref{criticalrlambda}). Fig. \ref{30degsint} illustrates this argument for $\alpha=\pi/3$. Shapes do exist for $\lambda<\lambda_{\text{min}}$ but without an inflection point. As $\lambda$ increases beyond $\lambda_{\text{min}}$, the inflection point moves from $r=1$ toward smaller $r$ and the slope at the inflection becomes increasingly vertical as $\lambda$ approaches $\lambda_c$. Fig. \ref{rmaxtheta} shows the location of the inflection (i.e., where $\sin\theta$ is maximum) versus $\lambda$.
\\
\\
The axial length of the rotating meniscus has two qualitatively different configuration types. The first type is associated to contact angles larger than zero. In these cases, the radial position of the maximum slope, $r_0$, corresponds to an inflection point, that is, $h'(r_0)>0, h''(r_0)=0$ (i.e., zero curvature) for all $\lambda<\lambda_c$; and the coefficient of the logarithmic divergence is $1/3$. In general, the location of the maximum slope changes with $\lambda$ for given contact angle $\alpha$, as shown for example in fig. \ref{30degsint}. In the $\alpha=\pi/2$ case analyzed in Sec. \ref{infiniteasympt}, however, the maximum slope location is independent of $\lambda$. The second type of behavior has a single element in the zero-contact angle case. This case always has infinite slope at $r_0=1$, but there is a non-zero curvature at the wall for all $\lambda<\lambda_c$, given by $r''(z)$, at the contact point $r=1$, see fig. \ref{shapesschem}. The coefficient of the logarithmic divergence for the zero-contact angle case is $2/\sqrt{3}$. Perhaps more significantly, this case is mathematically identical to that of the spinning bubble.\\
\\
Our results for the spinning bubble have distinctly practical implications for the measurement of surface tension by the spinning bubble tensiometer. In this method, the experimenter spins a container holding the gas-liquid pair of interest at high enough angular velocity $\omega$ that the bubble shape is nearly a straight cylinder. Then the radius, $R_b$, of the cylindrical bubble is measured and, assuming that $\lambda=4$, the corresponding surface tension, denoted $\sigma_4$ to reflect the assumption just made, is obtained from:
\beq
\lambda_c=4=\frac{\omega^2 R_b^3\rho}{\sigma_4}.
\label{errorsigma}
\eeq
However, for this to be valid the bubble must be close enough to its critical configuration which obtains when $4-\lambda\ll1$. Assessing this condition is not immediately obvious. In practice, an experimenter may perform several measurements at progressively higher angular velocities $\omega$, and evaluate $\sigma_4$ each time. The measurement would be satisfactory when the value of $\sigma_4$ attains the desired accuracy --for example by the first $n$ decimals remaining constant.\\
\\
While the method outlined above can be used to measure the surface tension, it does not inform about the value of $\lambda$ in a particular measurement. This is why the method is based on the \textit{assumption} that $\lambda=4$ and the experimenter must ensure that this condition is met with enough accuracy. The most direct way to find $\lambda$ is to evaluate $H$ from the ratio of bubble length to maximum radius, both of which can be measured. Using the theory, we may find $\lambda$ from the function $H(\lambda)$.\\
\\
Knowing the dimensional bubble volume, $\widetilde V$, from the difference of container volume and (incompressible) liquid volume, eq. (\ref{vasy}), rewritten here as
\beq
\widetilde V \sim R_b^3 (H\pi-4.18)
\label{vasy2}
\eeq
shows that, since $H$ diverges when $\lambda\to\lambda_c=4$ (eq. (\ref{Hasyspinning})), the bubble radius $R_b$ approaches zero as $\lambda\to4$. This is why the spinning bubble can only operate arbitrarily close to but not at $\lambda=4$. In practical terms, when an instrument spins to produce a bubble with, say, $H=10$ (i.e. the length is 10 times the maximum radius), the theory indicates that $\lambda\approx 4-2.8\times 10^{-3}$. Since $H(\lambda)$ is a measure of how close $\lambda$ is to 4, let us now consider the error in surface tension that one makes by assuming $\lambda=4$, as a function of $H$. This percent error, $\Delta$, is defined as
\beq
\Delta\equiv100\left(1-\frac{\sigma_4}{\sigma_{\text{actual}}}\right)=
     100\left(1-\frac{\lambda(H)}{4}\right).
\eeq
Joseph's \cite{joseph1994} argument that bubbles with $H>4$ can be considered to be at $\lambda=4$ is therefore inaccurate, as the present calculation predicts that $\lambda\approx3.48$ when $H=4$, yielding an error of 13\%. Thus, a bubble with $H=4$ is not long enough to be considered ``critical''. This is consistent with the plot of fig. \ref{VvsH} where $H=4$ is close to but not yet in the limiting long-$H$ regime where volume increases linearly with $H$. However, if we make the same assumption, $\lambda=4$, when the bubble length is $H=10$, the error drops to 0.07\%. The error of using eq. (\ref{errorsigma}) is given by $\Delta = 411.13 \,\exp(-\frac{\sqrt{3}}{2}H)$, with accuracy increasing with $H$.
\\
\\
In principle the experimenter need not assume $\lambda=4$, however, since we now have $H(\lambda)$ (i.e., the ratio of dimensional bubble length to maximum radius, both of which can be measured), described with a uniformly convergent asymptotic approximant (see Sec. \ref{approximants}) over the entire range $0\le\lambda<4$. Thus, the approximant for the bubble length  allows one to extend surface tension measurement to arbitrary values of $\lambda$ with just a simple evaluation of $H_A(\lambda)$ from measurements. In the absence of gravity, the approximant allows measurements in the intermediate-$\lambda$ region where the sensitivity to error in $H$ is still moderate. But working at a lower than critical $\lambda$ has the drawback that it would require precise measurements of both radius \emph{and} length.

\section{Summary and conclusion}
We have examined the problem of interface shapes in fluid systems under rigid-body rotation with a focus on finding exact solutions and the asymptotics of singular behaviors near $\lambda_c$. We studied two configurations of practical importance, e.g., a meniscus spanning the rotating container radius with arbitrary contact angle at the container wall; and a spinning bubble where the meniscus does not contact the container wall. Finding the asymptotic behavior of each meniscus configuration length as the critical rotation is approached as well as the series solution about $\lambda=0$ is important because such meniscus configurations arise in applications such as  a rotating reactor and a spinning bubble tensiometer. Knowing the form of the asymptotic divergence, one may construct efficient asymptotic approximants to evaluate each meniscus length at any rotation velocity uniformly and without solving a differential equation numerically.\\
\\
In conclusion, this work provides analyses that advance the interpretation of interface shapes of fluids in rigid body rotation. The analyses are strictily valid in zero-gravity, but their validity may be extended to normal gravity as long as the gravitational Bond number, $\rho g R^2/\sigma$, is much smaller than the rotational Bond number, $\lambda$. For two canonical configurations (meniscus spanning container radius and spinning bubble) we have found exact solutions for the axial meniscus length, $H$,  over the whole range of $\lambda$ and asymptotic behaviors near critical rotation. To remedy the poor convergence of the infinite-sum exact solutions we constructed convergent asymptotic approximants that greatly improve the convergence efficiency of the exact solution. Our results provide proof of concept of useful analytical calculation tools for applications ranging from controlling rotating reactors to measurement of surface tension with the spinning bubble method.

\section*{Acknowledgement}
Enrique Ram\'e is grateful to Dr. R. Balasubramaniam of Case Western Reserve University for help formulating initial ideas and calculations and for providing a critical sounding board during frequent discussions.

\begin{appendices}
\section{Explicit evaluation of $C_n$ in eq. (\ref{HTaylor2})}
To complement the recursive evaluation of $C_n$ of Sec. \ref{seriesbubble}, in this appendix we compute $C_n$ explicitly. Explicit expressions (as opposed to recursive relations) may be desirable for certain types of analysis.\\
\\
Because $f(r,\lambda)\equiv\sin(\theta)=r+\frac{\lambda}{8} (r-r^3)\le1$, we attempt expanding the denominator of $h'= f/\sqrt{1-f^2}$ as a series in powers of $f$. This approach is analogous to that employed in Sec. \ref{seriessolmeniscus} but now $\lambda$ is embedded in the square of the binomial $f$. Formally, this yields:
\beq
h'=\frac{f}{\sqrt{1-f^2}}=\sum_{m=0}^{\infty}f^{2m+1}\frac{\Gamma(\frac{1}{2}+m)}{\sqrt{\pi}\,\Gamma(m+1)}.
\label{hprimeser}
\eeq
Note that this is just the product of $f$ times the series of even powers for the denominator. Now expand the powers of $f$ using the binomial theorem:
\beq
f^{2m+1}=\sum_{i=0}^{2m+1}\frac{(2m+1)!}{i!\,(2m+1-i)!}
\left(\frac{\lambda}{8}\right)^{2m+1-i}r^i (r-r^3)^{2m+1-i}.
\eeq

A closed-form is available for the integration of the $r$-dependence in the above sum:
\beq
\int_0^1 f^{2m+1}\,dr=\sum_{i=0}^{2m+1}\frac{(2m+1)!}{i!\,(2m+1-i)!}
\left(\frac{\lambda}{8}\right)^{2m+1-i}
\frac{\Gamma(1+m)\Gamma(2m+2-i)}{2\,\Gamma(3m+3-i)}.
\eeq

We now use this to write the $r$-integral of eq. (\ref{hprimeser}) as:
\begin{align}
H(\lambda)&=2\sum_{m=0}^{\infty}\frac{\Gamma(\frac{1}{2}+m)\Gamma(2m+2)}{\sqrt{\pi}}
\sum_{i=0}^{2m+1}
\left(\frac{\lambda}{8}\right)^{2m+1-i}
\frac{1}{2\,\Gamma(i+1)\Gamma(3m+3-i)}\nonumber\\
&=\sum_{p=0}^{\infty} C_p \lambda^p.
\label{Hsum1}
\end{align}

It remains to extract the coefficient of $\lambda^p$, $C_p$. Let $2m+1-i=p$.  In the finite $i$-index sum, set $i=2m+1-p$ for each $m$. It follows that the coefficient of $\lambda^p$ is the result of an infinite sum:
\beq
C_p =\frac{1}{\sqrt{\pi}}\frac{1}{8^p}\sum_{m=m_{0}}^{\infty}\frac{\Gamma(\frac{1}{2}+m)\,\Gamma(2m+2)}{\Gamma(2m+2-p)\,\Gamma(2+m+p)}.
\label{cpofp}
\eeq
For a given $p$, the argument of $\Gamma(2m+2-p)$ in the denominator cannot be less than 1, i.e., $2m+1\ge p$. This sets the lowest $m$ in the sum, as $m\ge m_{0}\ge (p-1)/2$. When $p$ is odd, this condition sets the lowest $m$ directly; when $p$ is even, the lowest $m$ is the smallest integer that is larger than $(p-1)/2$.

\end{appendices}

\bibliography{Article.bib}
\end{document}